\documentclass[journal]{IEEEtran}
\IEEEoverridecommandlockouts
\usepackage{cite}
\usepackage{amsmath,amssymb,amsfonts}
\usepackage{algorithmic}
\usepackage{graphicx}
\usepackage{textcomp}
\usepackage{xcolor}

\def\BibTeX{{\rm B\kern-.05em{\sc i\kern-.025em b}\kern-.08em
    T\kern-.1667em\lower.7ex\hbox{E}\kern-.125emX}}

\begin{document}

\markboth{This paper has been accepted by IEEE ISGT NA 2019.}%
{Shell \MakeLowercase{\textit{et al.}}: Bare Demo of IEEEtran.cls for IEEE Journals}

\title{A Neural-Network-Based Optimal Control of\\ Ultra-Capacitors with System Uncertainties 
}

\author{\IEEEauthorblockN{Jiajun Duan\IEEEauthorrefmark{1},
		Zhehan Yi\IEEEauthorrefmark{1},
		Di Shi\IEEEauthorrefmark{1},
		Hao Xu\IEEEauthorrefmark{2}, and
		Zhiwei Wang\IEEEauthorrefmark{1}\\
		\IEEEauthorblockA{\IEEEauthorrefmark{1}GEIRI North America, San Jose, CA, 95134
		\\ Emails: \{jiajun.duan, zhehan.yi, di.shi, zhiwei.wang\}@geirina.net}\\
		\IEEEauthorblockA{\IEEEauthorrefmark{2}Department of Electrical and Biomedical Engineering\\
			University of Nevada,
			Reno, NA, 89557 \\} 
	\thanks{This work is supported by the SGCC Science and Technology Program \textit {Hybrid Energy Storage Management Platform for Integrated Energy System.}}		
}
}
\maketitle

\begin{abstract}
In this paper, a neural-network (NN)-based online optimal control method (NN-OPT) is proposed for ultra-capacitors (UCs) energy storage system (ESS) in hybrid AC/DC microgrids involving multiple distributed generations (e.g., Photovoltaic (PV) system, battery storage, diesel generator). Conventional control strategies usually produce large disturbances to buses during charging and discharging (C\&D) processes of UCs, which significantly degrades the power quality and system performance, especially under fast C\&D modes. Therefore, the optimal control theory is adopted to optimize the C\&D profile as well as to suppress the disturbances caused by UCs implementation. Specifically, an NN-based intelligent algorithm is developed to learn the optimal control policy for bidirectional-converter-interfaced UCs. The inaccuracies of system modeling are also considered in the control design. Since the designed NN-OPT method is decentralized that only requires the local measurements, plug \& play of UCs can be easily realized with minimal communication efforts.  In addition, the PV system is under the maximum power point tracking (MPPT) control to extract the maximum benefit. Both islanded and grid-tied modes are considered during the controller design. Extensive case studies have been conducted to evaluate the effectiveness of proposed method.
\end{abstract}


\section{Introduction}
Along with the development of smart grid, concerns from environmental protection, and considerations of future sustainability, the application of renewable energy resources have been greatly sped up in a form of microgrid \cite{01,di}. In the meanwhile, ESSs are being deployed to compensate for the intermittency of renewable DGs such as PV and wind. The wide employment of ESS can enable a better and more reliable power supply, avoid the extraordinary updating cost of the conventional power system, and provide an emergency backup after a catastrophic incident. Additionally, it may entirely alternate the traditional supply and demand patterns \cite{02}.\par

From this perspective, UC, which has been widely applied in microgrids, is regarded as one of the most promising ESS due to its flexible size and high power density \cite{03}. However, the fast dynamics have posted several challenges to the smooth C\&D control of UCs \cite{04}. In general, a UC can be considered as a generator during the discharging mode, and as a load during the charging mode. For the rest of time, UC should be isolated from the system to prevent the continuous repetition of C\&D \cite{04}. However, the initiating or switching of C\&D modes in microgrids may result in considerable disturbances on the common bus voltage, especially under the fast C\&D mode \cite{Sheng1}. It is harmful to the sensitive loads (e.g., data center) on the bus and may further trigger false protection actions \cite{05,Yi1}. This is opposite to the original design purpose of ESS, i.e., to decrease the system disturbance and improve the power quality. Therefore, a desired control policy of UC should provide a smooth C\&D solution to power systems as well as possess plug \& play capability for DGs and loads \cite{Sheng2}.\par

At present, constant charging current control remains as one of the most popular methods for ESSs \cite{12, 06}. After initializing, the UC is charged or discharged at the maximum constant rate. Then, it will be turned off forcibly once the desired state of charge (SOC) is reached. The major challenge of the constant charging current control is that it will produce large disturbances at both beginning and endding stages of the C\&D for UC. Another widely employed method is the proportion-integration (PI)-based control \cite{07,Yi2}, which also has a quite high industrial maturity. Conventional PI-based control methods mainly have two control loops, i.e. the inner voltage and the outer current loops. However, in order to avoid the over C\&D problem, the inner voltage loop is usually designed as a P controller. This greatly degrades the overall performance of system and may introduce a large disturbance during the C\&D processes, which has been thoroughly illustrated in previous work \cite{04} and later in this one. \par

In order to address the aforementioned challenges, a real-time NN-OPT control algorithm is proposed to optimize the C\&D process of UC inside a hybrid AC/DC microgrid. Due to the existence of nonlinearity and inaccuracies during the dynamic system modeling process, an NN is developed to learn the optimal control input through online tuning. The UC is connected with a PV array and diesel generator (DG) through a modified bidirectional power converter (BPC). During control design, both islanded and grid-tied scenarios are considered. In islanded mode, the diesel generator is responsible for maintaining the bus voltage at the point of common coupling  (PCC), while the power exchange between AC and DC network is controlled by a voltage source converter (VSC). In grid-tied mode, bus voltage at PCC is regulated by the main grid and reactive power is controlled by the VSC. In either scenario, PV generation is working in the maximum power point tracking (MPPT) mode using incremental conductance and integral regulator (IC\&IR) technique. The effectiveness of proposed control algorithm is demonstrated through extensive simulations and the overall performance is proven to have a significant improvement comparing to the conventional methods. \par

The rest of paper is organized as follows. In Section II, the problem formulation for the studied system is introduced. Then,  the proposed learning based control algorithm and the detailed implementation procedures are given in Section III.  The case study results and corresponding analysis are illustrated in Section IV. Finally, the conclusion and recommended future work are presented in Section V. \par

 \section{Problem Formulation}

\begin{figure}
	\centering
	\includegraphics[width=\linewidth]{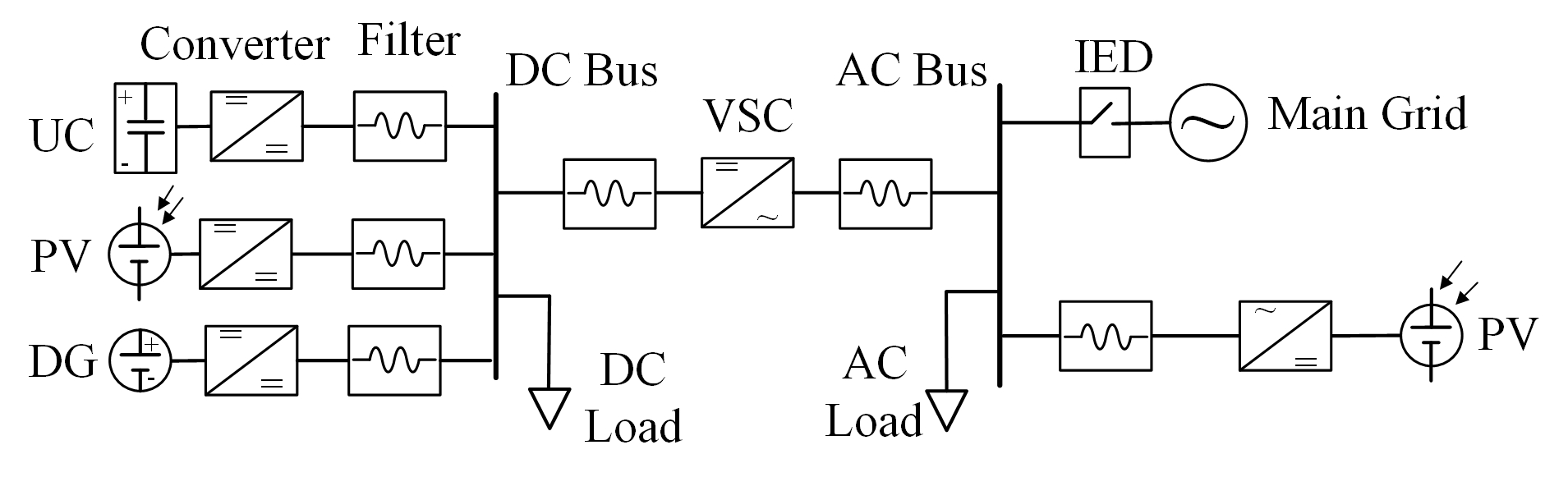}
	\caption{Diagram of the considered hybrid AC/DC microgrid.}
	\label{fig1}
\end{figure}

\begin{figure}
	\centering
	\includegraphics[width=0.85\columnwidth]{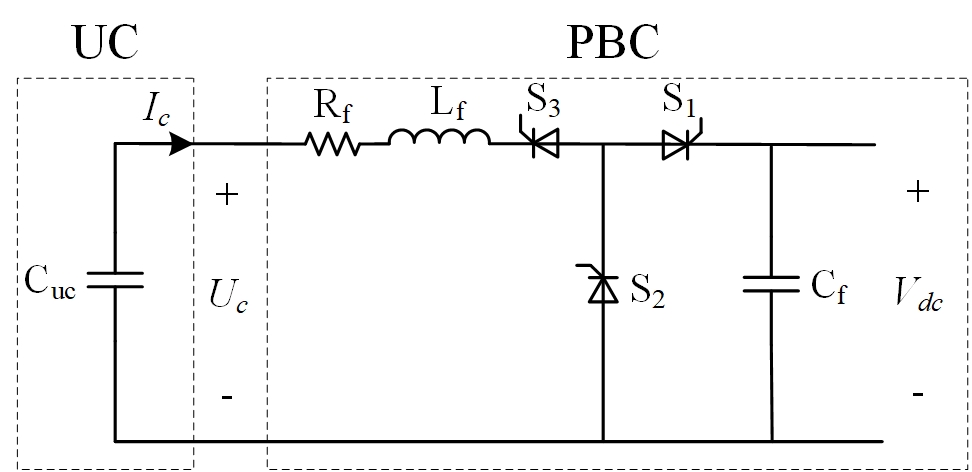}
	\caption{Switch-level UC interface circuit.}
	\label{fig2}
\end{figure}

As is presented in Fig. \ref{fig1}, a microgrid consisting of both AC and DC buses is connected to the main grid through an intelligence electric device (IED). The PV array is operating in MPPT mode with a UC as the energy storage. The switch-level circuit of UC is illustrated in Fig. \ref{fig2}. The UC is connected to the DC PCC through a modified BPC. The switch $\rm{S_{3}}$ is added to prevented self-discharging of UC during the initializing period of C\&D. In the normal C\&D process, if $\rm{S_{1}}$ is off and $\rm{S_{2}}$ is under control, the PBC is working as a boost converter to discharge the UC. If $\rm{S_{2}}$ is off and $\rm{S_{1}}$ is under control, the PBC is working as a buck converter to charge the UC.\par

In general, the simplified dynamics of a UC can be represented by Eqn. (\ref{eq1})
\begin{equation}
\label{eq1}
\dot{U}_{c}(t)=\pm\frac{1}{C_{uc}}I_{c}(t)
\end{equation}
where $U_{c}$ is the UC terminal voltage; $I_{c}$ is the charging (-) or discharging (+) current; $C_{uc}$ is the UC capacitance. Since Eqn. (\ref{eq1}) only represents an ideal model, the error dynamic of Eqn. (\ref{eq1}) considering the unmodelled part of system can be expressed as
\begin{equation}
\label{eq2}
\dot{e}(t)=g(e(t))u(t)+D(t), e(0)=e_{0}
\end{equation}
where  $e(t)=U_{c}(t)-U_{c}^{*}$ with $U_{c}^{*}$ being SOC reference, $u(t)=I_{c}(t)$ being the control input, $D(t)$ being the non-linear perturbation of dynamic system bounded by $\|D(t)\|\leq{d_{max}}$.\par
To reduce the disturbances caused by C\&D, the optimal control problem will be formulated with respect to $e(t)$ and $u(t)$. In the next section, an online learning-based optimal control method for the uncertain systems is developed.\par

\section{NN-Based Optimal Control Design for UC}
In this section, the optimal control for the uncertain nonlinear system (\ref{eq2}) is firstly formulated. Then, an NN-based algorithm is applied to learn the optimal control input in an online training fashion. \par
\subsection{Optimal Control Formulation for Uncertain System}
For a nonlinear system without uncertainty $D(t)$, i.e., $\dot{e}(t)=g(e(t))u(t)$, the infinite-horizon integral cost function can be designed as
\begin{equation}
\label{eq3}
J(e_{0}, u)=\int_{0}^{\infty}r(e, u)dt
\end{equation}
where $r(e, u)=Q(e)+u^{T}Ru$, $Q(e)$ is a positive definite function of $e$ and $R$ is a symmetric positive definite matrix. According to \textit{Theorem 1} in \cite{08}, control law $u(e)$ can guarantee the asymptotic stability of a closed-loop system if there exist a positive definite continuously differentiable function $V(e)$, a bounded function $\Gamma(e)$ and a feedback control law $u(e)$ such that 
\begin{equation}
\label{eq4}
\begin{cases}
V_{\partial e}^{T}D(t)\le\Gamma(e)\\
V_{\partial e}^{T}[g(e)u]+\Gamma(e)+Q(e)+u^{T}Ru=0
\end{cases}
\end{equation}
where $V_{\partial e}$ is the partial derivative of the cost function $V(e)$ with respect to $e$. Then, cost function Eqn. (\ref{eq3}) satisfies
\begin{equation}
\label{eq5}
\sup_{D(t)\in M}J(e_{0}, u)\le J_{d}(e_{0}, u)=V(e_{0})
\end{equation}
where ``sup'' denotes the supremum operator that finds the minimal cost $J_{d}(e_{0}, u) \ge J(e_{0}, u)$ for any perturbation $D(t)\in M, M=\{ D(t)|D(t)\in\Re, \|D(t)\|\leq{d_{max}} \}$, where
\begin{equation}
\label{eq6}
J_{d}(e_{0}, u)=\int_{0}^{\infty}[r(e, u)+\Gamma(e)]dt
\end{equation}
represents the modified cost function considering the system uncertainties. Accordingly, the bounded function $\Gamma(e)$ can be designed as
\begin{equation}
\label{eq7}
\Gamma(e)=\frac{1}{4}V_{\partial e}^{T}V_{\partial e}+d_{max}^{2}
\end{equation}

Based on the definition of $\Gamma(e)$, cost function Eqn. (\ref{eq6}) can be rewritten as
\begin{eqnarray}
\label{eq8}
\begin{split}
J_{d}&(e_{0}, u)=V(e_{0})\\
&=\int_{0}^{T}[r(e, u)+\Gamma(e)]dt+\int_{T}^{\infty}[r(e, u)+\Gamma(e)]dt\\
&=\int_{0}^{T}[r(e, u)+\Gamma(e)]dt+V(e)
\end{split}
\end{eqnarray}

\noindent Since $V(e)$ is continuously differentiable, Eqn. (\ref{eq8}) becomes
\begin{eqnarray}
\label{eq9}
\begin{split}
\lim_{T\to0}&\frac{V(e_{0})-V(e)}{T}=\lim_{T\to0}\frac{1}{T}\int_{0}^{T}[r(e, u)+\Gamma(e)]dt\\
&\Rightarrow\dot{V}(e)=V_{\partial e}^{T}[g(e)u+D]=-r(e, u)-\Gamma(e)\\
&\Rightarrow 0=V_{\partial e}^{T}[g(e)u+D]+r(e, u)+\Gamma(e)
\end{split}
\end{eqnarray}
which is an infinitesimal version of Eqn. (\ref{eq8}) as well as a non-linear Lyapunov equation.\par
Based on Eqn. (\ref{eq9}), Hamiltonian of the optimal control problem can be defined as
\begin{equation}
\label{eq10}
H(e, u, V_{\partial e})=Q(e) +u^{T}Ru+V_{\partial e}^{T}[g(e)u+D]+\Gamma(e)
\end{equation}

\noindent The optimal cost function $V^{*}(e)=\min_{u\in\Omega}\int_{0}^{T}[r(e, u)+\Gamma(e)]dt$ can be solved as the solution of Hamilton-Jacobi-Bellman (HJB) equation
\begin{equation}
\label{eq11}
\min_{u\in\Omega}H(e, u, V_{\partial e}^{*})=0
\end{equation}

\noindent By solving $\partial H(e, u, V_{\partial e}^{*})/\partial u=0$, the optimal control law $u^{*}$ can be derived as
\begin{equation}
\label{eq12}
u^{*}=-\frac{1}{2}R^{-1}g^{T}(e)V_{\partial e}^{*}
\end{equation}

\noindent Substituting Eqn. (\ref{eq12}) into Eqn. (\ref{eq11}), the HJB equation in terms of $V_{\partial e}^{*}$ can be represented as
\begin{eqnarray}
\label{eq13}
\begin{split}
0=~&Q(e) +\frac{1}{4}V_{\partial e}^{*T}V_{\partial e}^{*}+d_{max}^{2}+V_{\partial e}^{*T}D\\
&-\frac{1}{4}V_{\partial e}^{*T}g(e)R^{-1}g(e)^{T}V_{\partial e}^{*}
\end{split}
\end{eqnarray}

\textit{Theorem 1 (Optimal Control Law $u^{*}$)}: For a non-linear uncertain system given in Eqn. (\ref{eq2}) with the cost function defined in Eqn. (\ref{eq5}) and HJB equation defined in Eqn. (\ref{eq13}), provided any admissible control $u$, the cost function Eqn. (\ref{eq5}) is smaller than a guaranteed cost bound $J_{b}$ given as
\begin{equation}
\label{eq14}
J_{b}=V^{*}(e_{0})+\int_{0}^{T}(u-u^{*})^{T}R(u-u^{*})dt
\end{equation}
where $u^{*}$ is obtained from Eqn. (\ref{eq12}). If $u=u^{*}$, the cost $J_{b}$ is guaranteed to be minimized, i.e., $J_{b}=V^{*}(e_{0})$.\par
\textit{Proof}: According to Eqn. (\ref{eq9}) and the definition of $V^{*}(e)$, the cost function Eqn. (\ref{eq5}) with respect to any arbitrary $u$ can be rewritten as
\begin{equation}
\label{eq15}
J(e_{0}, u)=V^{*}(e_{0})+\int_{0}^{T}[r(e,u)+\dot{V}^{*}(e)]dt
\end{equation}

\noindent Using Eqn. (\ref{eq9}) and Eqn. (\ref{eq13}), one can obtain that
\begin{eqnarray}
\label{eq16}
\begin{split}
r(e,u)&+\dot{V}^{*}(e)=Q(e) +u^{T}Ru+V_{\partial e}^{*T}[g(e)u+D]\\
=&u^{T}Ru+V_{\partial e}^{*T}g(e)u+\frac{1}{4}V_{\partial e}^{*T}g(e)R^{-1}g(e)^{T}V_{\partial e}^{*}\\
&-\frac{1}{4}V_{\partial e}^{*T}V_{\partial e}^{*}-d_{max}^{2}\\
\le & u^{T}Ru+V_{\partial e}^{*T}g(e)u+\frac{1}{4}V_{\partial e}^{*T}g(e)R^{-1}g(e)^{T}V_{\partial e}^{*}
\end{split}
\end{eqnarray}

Recalling Eqn. (\ref{eq12}), Eqn. (\ref{eq16}) can be compiled into a square form with respect to $R^{-1}g(e)^{T}V_{\partial e}^{*}/2$ as  
\begin{equation}
\label{eq17}
r(e,u)+\dot{V}^{*}(e)\le(u-u^{*})^{T}R(u-u^{*}) 
\end{equation}
which implies that Eqn. (\ref{eq14}) holds. Thus, if $u=u^{*}$, the cost $J_{b}$ is guaranteed to be minimized, i.e., $J_{b}=V^{*}(e_{0})$, and the corresponding optimal control input is $u^{*}$. Proof completed.$\Diamond$ \par

Obviously, the solution for HJB Eqn. (13) is required to obtain the optimal control input $u^{*}$. However, it is quite challenging to solve a non-linear partial derivative function. Therefore, an NN-based adaptive method is developed in the next subsection to approximate the solution of HJB equation.\par
\subsection{NN-based Optimal Control Design}
According to the universal approximation property of NN \cite{09}, cost function $V^{*}(e)$ can be expressed by an NN with the activation function on a compact set $\Omega$ in form as
\begin{equation}
\label{eq18}
V^{*}(e)=W_{c}^{*T}\sigma(e, u)+\epsilon
\end{equation}
where $W_{c}\in\Re^{n}$ is the target NN weight, $\sigma(e, u)\in\Re^{n}$ is the bounded activation function containing $n$ neurons, and $\epsilon$ is the NN reconstruction error. The target NN weights $W_{c}$ and reconstruction error $\sigma(e, u)$ are assumed to be bounded by $\|W_{c}\|\le W_{M}$ and $\|\epsilon\|\le \epsilon_{M}$ \cite{09}. \par
Since the target NN weight $W_{c}$ is unknown, an NN-based estimator is designed to approximate the cost function as
\begin{equation}
\label{eq19}
\hat{V}(e)=\hat{W}_{c}^{T}\sigma(e, u)
\end{equation}
where $\hat{V}(e)$ is the approximated cost function and $\hat{W}_{c}\in\Re^{n}$ is the estimated NN weight.
According to Eqn. (\ref{eq12}) and (\ref{eq19}), the estimate optimal control is given by
\begin{equation}
\label{eq20}
\hat{u}=-\frac{1}{2}R^{-1}g^{T}(e)\triangledown\sigma(e, u)^{T}\hat{W}_{c}
\end{equation}
where $\triangledown\sigma(e, u)=\partial \sigma(e, u)/\partial e$. Substituting Eqn. (\ref{eq19}) and (\ref{eq20}) into (\ref{eq10}), the approximated Hamiltonian can be obtained as
\begin{eqnarray}
\label{eq21}
\begin{split}
\hat{H}(e, \hat{u}, \hat{W}_{c})=&Q(e) -\frac{1}{4}\hat{W}_{c}^{T}\triangledown\sigma g(e)R^{-1}g(e)^{T}\triangledown\sigma^{T}\hat{W}_{c}\\
&+\frac{1}{4}\hat{W}_{c}^{T}\triangledown\sigma \triangledown\sigma^{T}\hat{W}_{c}+d_{max}^{2}
\end{split}
\end{eqnarray}

Due to the system uncertainty and cost function estimation error, the estimated Hamiltonian cannot hold, i.e., $\hat{H}(e, \hat{u}, \hat{W}_{c})\neq0$. According to the optimal control theory \cite{10}, the estimated cost function can converge close to the ideal target while the approximated Hamiltonian equation approaches to the ideal Hamiltonian,
i.e., $\hat{H}(e, \hat{u}, \hat{W}_{c})\to H(e, u^{*}, W_{c}^{*})=0$. Inspired by this, the updating law for tuning the NN weight of cost function estimator can be designed as
\begin{eqnarray}
\label{eq22}
\begin{split}
\dot{\hat{W}}_{c}=\frac{\alpha_{1}}{2}\Theta(e, \hat{u})\triangledown\sigma g(e)R^{-1}g(e)^{T}J_{1\partial e}-\frac{\alpha_{2}\omega\hat{H}}{(1+\omega^{T}\omega)^{2}}
\end{split}
\end{eqnarray}
where $\alpha_{1}>0$ and $\alpha_{2}>0$ are designed control coefficients, $\omega=-[\triangledown\sigma g(e)R^{-1}g(e)^{T}\triangledown\sigma^{T}\hat{W_{c}}]/2$, and $\Theta(e, \hat{u})$ is an index operator given by
\begin{equation}
\label{eq23}
\Theta(e, \hat{u})=
\begin{cases}
0, ~\forall \dot{J}_{1}=J_{1\partial e}^{T}\dot{e}<0\\
1, ~\text{otherwise}
\end{cases}
\end{equation}
where $J_{1}$ is a unbounded Lyapunov candidate and $J_{1\partial e}$ is its partial derivative with respect to $e$. Moreover,  $J_{1\partial e}$ can be defined similar to \cite{11} as
\begin{equation}
\label{eq24}
\|\dot{e}\| \le c_{1}\|e\|\equiv (c_{2}\|J_{1\partial e}\|)^{\frac{1}{4}}
\end{equation}

where $c_{1}$ and $c_{2}$ are constants. Note that $\|J_{1\partial e}\|$ can be selected to satisfy the general bound, e.g., $J_{1}=\frac{1}{5}(e^{T}e)^{\frac{5}{2}}.$\par

\textit{Theorem 2 (Convergence of the Optimal Control)}: Consider the nonlinear uncertain system in Eqn. (\ref{eq2}) with control law in Eqn. (\ref{eq20}) and NN weight updating law in Eqn. (\ref{eq22}), both the tracking error $e$ and the weight estimation error of NN (i.e. $\tilde{W}_c=W_{c}^{*}-\hat{W}_{c})$  are guaranteed to be uniformly ultimately bounded.\par
\textit{Proof}: Omitted in here and can be refereed to \cite{08,11}.\par
\begin{figure}
	\centering
	\includegraphics[width=\linewidth]{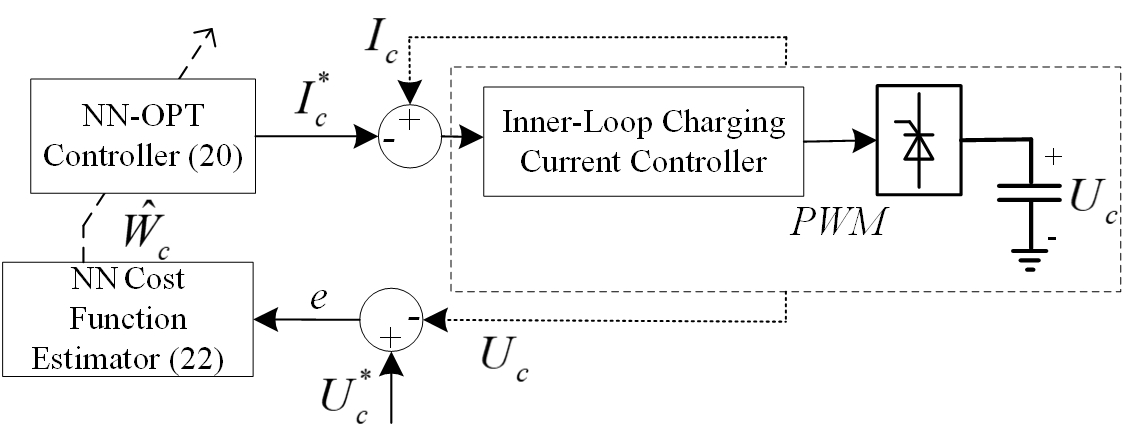}
	\caption{Schematic of the control implementation.}
	\label{fig3}
\end{figure}

 
The implementation of NN-based optimal controller is presented in Fig. \ref{fig3}. The proposed controller scales the voltage measurement of UC as the SOC and calculates the control input. Then, the inner loop current controller is responsible for tracking the optimal charging current reference using PI or hysteresis-based algorithm. The rest of system including DG and VSC are under the traditional PI-based controllers with control objectives described in Section I.  \par

\section{Case Studies}
In order to evaluate the effectiveness of proposed method, in this section, the NN-OPT controller is implemented in the system as described in Fig. \ref{fig1} with a sampling frequency of 10\,kHz. The simulation is conducted with a detailed switching-level system modeled in the Matlab/Simulink Simpowersystem toolbox. The controller and system parameters are given in Table \ref{table1}. Both islanded and grid-tied modes are presented in this section. Benchmarking studies are carried out by comparing the proposed method with the conventional PI-based method presented in [4]. \par
\vspace{0.5em}
\begin{table}
	\centering 
	\renewcommand{\arraystretch}{1.3}
	\caption{Case Study System Parameters}
	\label{table1}
	\begin{tabular}{p{1.7cm}p{1.7cm}p{1.7cm}p{1.7cm}}
		\hline \hline 
		Parameter& Value & Parameter & Value \\ 
		\hline 
		$C_{uc}$ & 5.7 F& $C_{f}$&20 $\mu$F \\ 
		$R_{f}$ & 0.1 $\Omega$& $L_{f}$&2 mH \\ 
		$\alpha_{1}$ & 2 & $\alpha_{2}$ & 0.24 \\ 
		\hline	\hline 
	\end{tabular} 
\end{table}
\subsection{Case I. Islanded Mode}
In this case, the proposed NN-OPT control method is firstly tested in an islanded microgrid. The SOC reference of UC, which can be represented by $U_{c}^{*}$ \cite{05}, is set to be discharged from 30V to 29V. The DC bus voltage is maintained at 48V constantly by a diesel generator. \par 
\begin{figure}
	\centering
	\includegraphics[width=\linewidth]{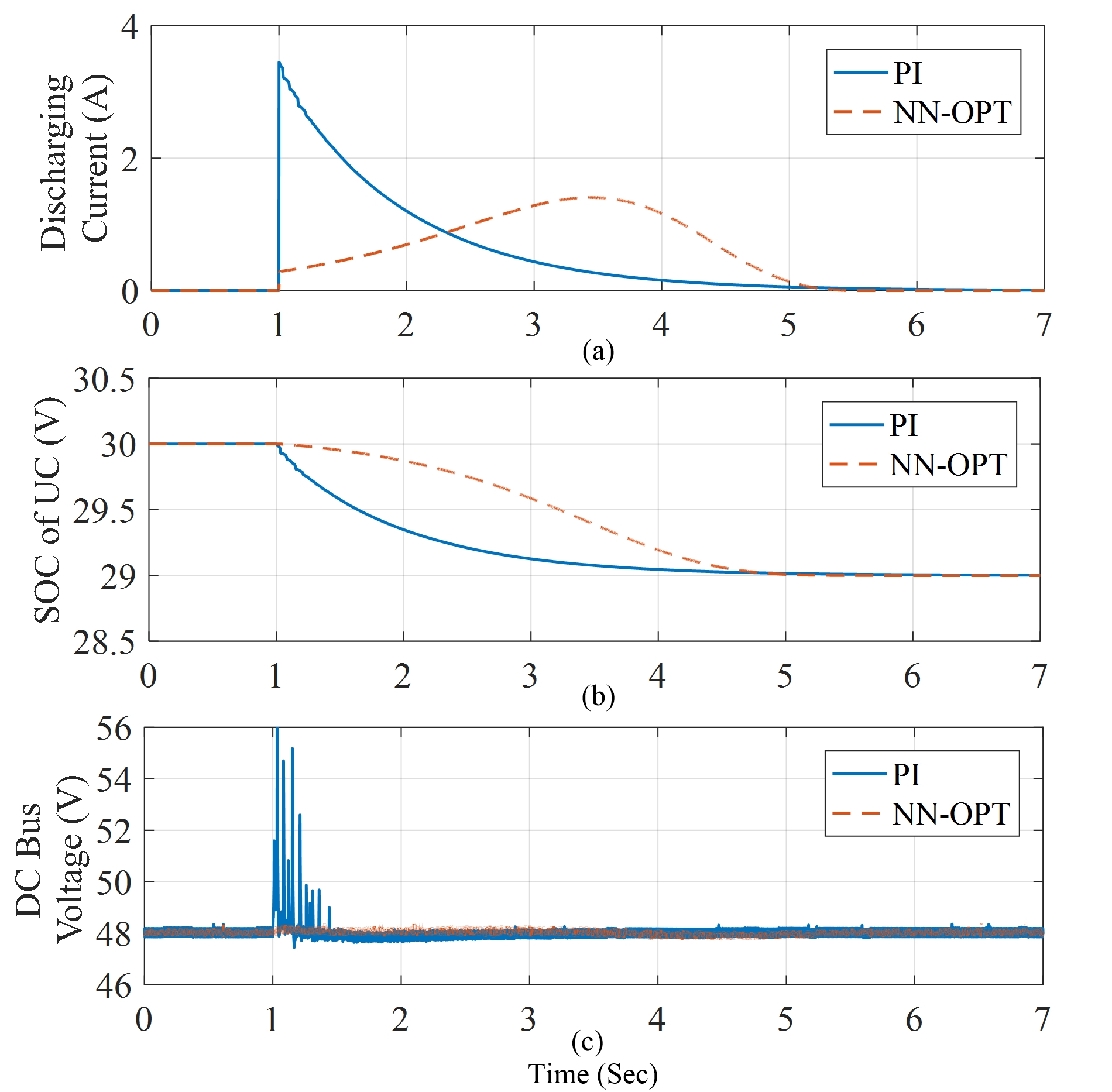}
	\caption{Simulation results in islanded mode (\itshape{Case I}).}
	\label{fig4}
\end{figure}
The results are presented in Fig. \ref{fig4}, among which Fig. 4(a) is the discharging current, Fig. 4(b) is the SOC of UC, and Fig. 4(c) is the voltage response of the DC bus. As can be seen, the conventional PI-based controller introduces an abrupt discharging current when UC is initialized to work, which causes an unexpected disturbance to the system. Consequently, the DC bus voltage appears multiple spikes over 8V (about 17\%). This is harmful to sensitive loads connected to the DC bus and is large enough to trigger the false protection action under certain circumstances \cite{05}. On the contrary, the proposed optimal control greatly eliminates the overshoot of discharging current during the initializing period. It also directly optimizes the entire C\&D process of UC which is reflected by the SOC profile. As a result, except for the normal harmonic caused by switching devices, there is barely any disturbances that can be observed on the DC bus voltage.\par

\begin{figure}
	\centering
	\includegraphics[width=\linewidth]{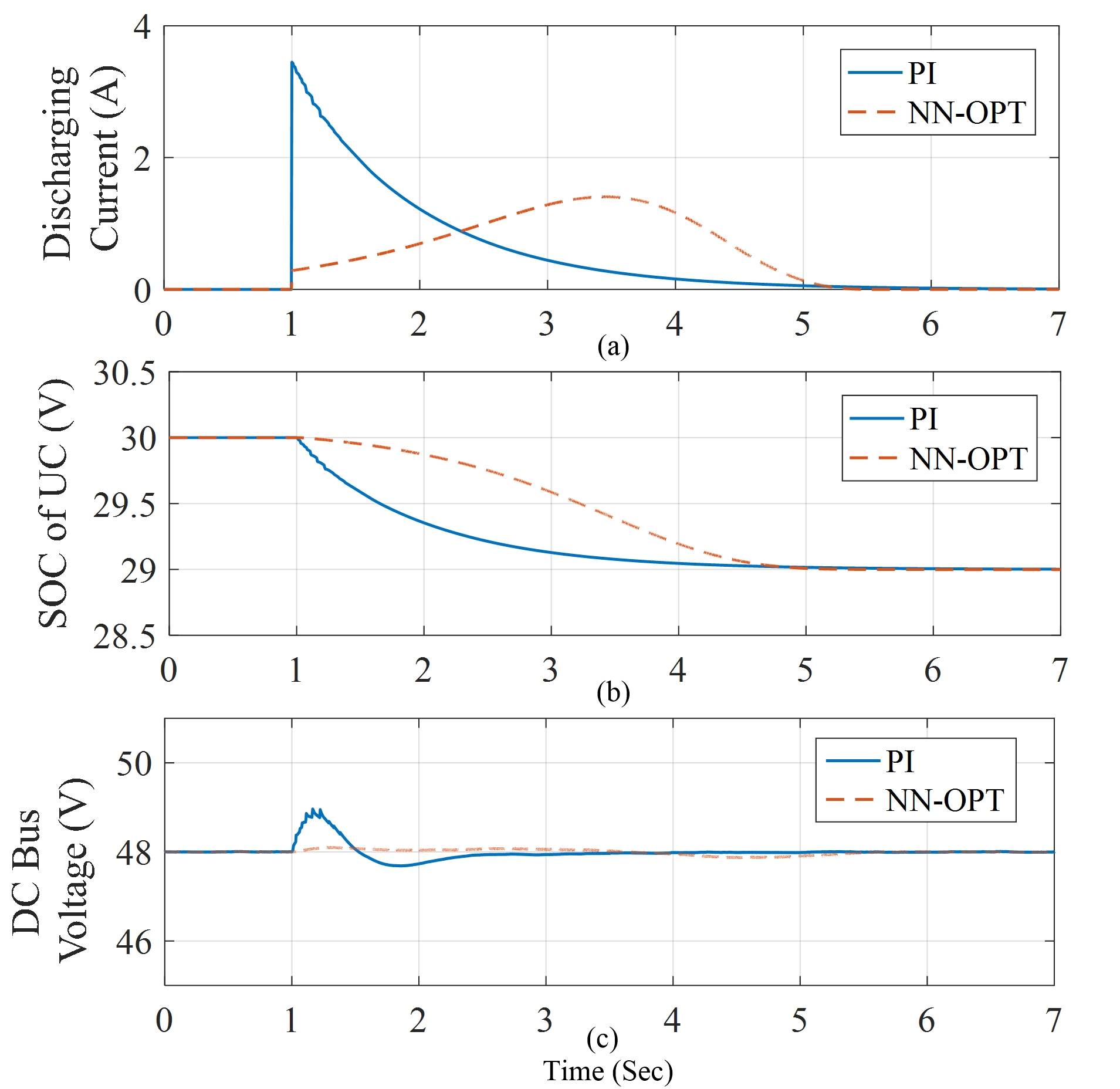}
	\caption{Simulation results in grid-tied mode (\itshape{Case II}).}
	\label{fig5}
\end{figure}

\subsection{Case II. Grid-tied Mode}
The proposed control method is also tested in the grid-tied situation under the same settings as \textit{Case I}. This time, the DC bus voltage and reactive power are regulated by the main grid through a VSC. The responses of discharging current, SOC of UC, and DC bus voltage are shown in Fig. 5(a), Fig. 5(b) and Fig. 5(c), respectively. As can be observed, in grid-tied mode, the main grid can provide a stronger voltage support comparing to the microgrid in islanded mode. Therefore, under the same discharging profile, both harmonics and oscillations on DC bus voltage are significantly reduced. Since the proposed control method has reached the optimal point, there is not much improvement can be observed. Even though, the proposed NN-OPT method exhibits a superior performance to conventional PI-based controller, which is directly reflects by the DC bus voltage (in Fig. 5(c)). Accordingly, the effectiveness of presented method is fully demonstrated. \par

\section{Conclusion}
In this paper, a novel NN-OPT control method is proposed for UC storage system in hybrid AC/DC microgrids involving PV and diesel generator. Firstly, the UC is modeled as a nonlinear system with uncertainties. Then, an NN is developed to learn the online optimal control input for such system. Utilizing optimal control theory, the C\&D profile of UC is optimized to suppress the disturbances caused by integration of ESS. The proposed method is tested in a hybrid microgrid under both islanded and grid-tied modes. The results demonstrate that the developed control method has a significant improvement comparing to the conventional controller. It is noteworthy that the proposed method is a unified approach that can be scaled to various sizes and types of ESS (e.g., Li-ion batteries, hybrid ESS). Therefore,  in the future, this work can be promoted to a more general ESS control system. In addition, the hardware experimentation can be conducted to further evaluate effectiveness of the proposed scheme. \par


\begin{thebibliography}{10}
	\providecommand{\url}[1]{#1}
	\csname url@samestyle\endcsname
	\providecommand{\newblock}{\relax}
	\providecommand{\bibinfo}[2]{#2}
	\providecommand{\BIBentrySTDinterwordspacing}{\spaceskip=0pt\relax}
	\providecommand{\BIBentryALTinterwordstretchfactor}{4}
	\providecommand{\BIBentryALTinterwordspacing}{\spaceskip=\fontdimen2\font plus
		\BIBentryALTinterwordstretchfactor\fontdimen3\font minus
		\fontdimen4\font\relax}
	\providecommand{\BIBforeignlanguage}[2]{{%
			\expandafter\ifx\csname l@#1\endcsname\relax
			\typeout{** WARNING: IEEEtran.bst: No hyphenation pattern has been}%
			\typeout{** loaded for the language `#1'. Using the pattern for}%
			\typeout{** the default language instead.}%
			\else
			\language=\csname l@#1\endcsname
			\fi
			#2}}
	\providecommand{\BIBdecl}{\relax}
	\BIBdecl
	
	\bibitem{01}
	J.~Duan, C.~Wang, H.~Xu, and W.~Liu, ``Distributed control of
	inverter-interfaced microgrids based on consensus algorithm with improved
	transient performance,'' \emph{IEEE Transactions on Smart Grid}, pp. 1--1,
	2018.
	
	\bibitem{di}
	D.~Shi, X.~Chen, Z.~Wang, X.~Zhang, Z.~Yu, X.~Wang, and D.~Bian, ``A
	distributed cooperative control framework for synchronized reconnection of a
	multi-bus microgrid,'' \emph{IEEE Transactions on Smart Grid}, pp. 1--1,
	2018.
	
	\bibitem{02}
	V.~A. Boicea, ``Energy storage technologies: The past and the present,''
	\emph{Proceedings of the IEEE}, vol. 102, no.~11, pp. 1777--1794, Nov 2014.
	
	\bibitem{03}
	P.~J. Grbovic, P.~Delarue, P.~L. Moigne, and P.~Bartholomeus, ``Modeling and
	control of the ultracapacitor-based regenerative controlled electric
	drives,'' \emph{IEEE Transactions on Industrial Electronics}, vol.~58, no.~8,
	pp. 3471--3484, Aug 2011.
	
	\bibitem{04}
	W.~Im, C.~Wang, L.~Tan, and W.~Liu, ``Cooperative controls for pulsed power
	load accommodation in a shipboard power system,'' \emph{IEEE Transactions on
		Power Systems}, vol.~31, no.~6, pp. 5181--5189, Nov 2016.
	
	\bibitem{Sheng1}
	Y.~Luo, S.~Srivastava, M.~Andrus, and D.~Cartes, ``Application of distubance
	metrics for reducing impacts of energy storage charging in an mvdc based
	ips,'' in \emph{2013 IEEE Electric Ship Technologies Symposium (ESTS)}, April
	2013, pp. 287--291.
	
	\bibitem{05}
	J.~Duan, C.~Wang, and H.~Xu, ``Distributed control of inverter-interfaced
	microgrids with bounded transient line currents,'' \emph{IEEE Transactions on
		Industrial Informatics}, vol.~14, no.~5, pp. 2052--2061, May 2018.
	
	\bibitem{Yi1}
	Z.~Yi and A.~H. Etemadi, ``Line-to-line fault detection for photovoltaic arrays
	based on multiresolution signal decomposition and two-stage support vector
	machine,'' \emph{IEEE Transactions on Industrial Electronics}, vol.~64,
	no.~11, pp. 8546--8556, Nov 2017.
	
	\bibitem{Sheng2}
	Y.~Luo, C.~Wang, L.~Tan, G.~Liao, M.~Zhou, D.~Cartes, and W.~Liu, ``Application
	of generalized predictive control for charging super capacitors in microgrid
	power systems under input constraints,'' in \emph{2015 IEEE International
		Conference on Cyber Technology in Automation, Control, and Intelligent
		Systems (CYBER)}, June 2015, pp. 1708--1713.
	
	\bibitem{12}
	V.~Vu, D.~Tran, and W.~Choi, ``Implementation of the constant current and
	constant voltage charge of inductive power transfer systems with the
	double-sidedlcccompensation topology for electric vehicle battery charge
	applications,'' \emph{IEEE Transactions on Power Electronics}, vol.~33,
	no.~9, pp. 7398--7410, Sept 2018.
	
	\bibitem{06}
	X.~Zheng, X.~Liu, Y.~He, and G.~Zeng, ``Active vehicle battery equalization
	scheme in the condition of constant-voltage/current charging and
	discharging,'' \emph{IEEE Transactions on Vehicular Technology}, vol.~66,
	no.~5, pp. 3714--3723, May 2017.
	
	\bibitem{07}
	X.~Feng, J.~Hu, Y.~Tao, H.~Liu, and D.~Liu, ``Research of off-grid energy
	storage converter based on repetitive control and pi control,'' in \emph{2016
		China International Conference on Electricity Distribution (CICED)}, Aug
	2016, pp. 1--4.
	
	\bibitem{Yi2}
	Z.~Yi, W.~Dong, and A.~H. Etemadi, ``A unified control and power management
	scheme for {PV}-battery-based hybrid microgrids for both grid-connected and
	islanded modes,'' \emph{IEEE Transactions on Smart Grid}, vol.~PP, no.~99,
	pp. 1--1, 2017.
	
	\bibitem{08}
	Y.~Huang, ``Optimal guaranteed cost control of uncertain non-linear systems
	using adaptive dynamic programming with concurrent learning,'' \emph{IET
		Control Theory Applications}, vol.~12, no.~8, pp. 1025--1035, 2018.
	
	\bibitem{09}
	F.~Lewis, S.~Jagannathan, and A.~Yesildirak, \emph{Neural network control of
		robot manipulators and non-linear systems}.\hskip 1em plus 0.5em minus
	0.4em\relax CRC Press, 1998.
	
	\bibitem{10}
	F.~L. Lewis, D.~Vrabie, and V.~L. Syrmos, \emph{Optimal control}.\hskip 1em
	plus 0.5em minus 0.4em\relax John Wiley \& Sons, 2012.
	
	\bibitem{11}
	T.~Dierks and S.~Jagannathan, ``Optimal control of affine nonlinear
	continuous-time systems,'' in \emph{American Control Conference (ACC),
		2010}.\hskip 1em plus 0.5em minus 0.4em\relax IEEE, 2010, pp. 1568--1573.
	
\end{thebibliography}

\end{document}